\documentclass[a4paper,11pt,fleqn]{article}

\usepackage{amssymb}

\usepackage{lmodern}

\usepackage[T1]{fontenc}

\usepackage{fullpage}
\usepackage{graphicx}
\usepackage{amssymb}
\usepackage{amsfonts}
\usepackage{amsthm}
\usepackage{amsmath}
\usepackage{cleveref}

\swapnumbers

\newtheorem{definition}{Definition}[section]

\newtheorem{theo}[definition]{Theorem}
\newtheorem{lemma}[definition]{Lemma}

\newtheorem{prop}[definition]{Proposition}

\DeclareMathOperator{\RE}{Re}
\DeclareMathOperator{\IM}{Im}

\DeclareMathOperator{\tr}{tr}
\DeclareMathOperator{\ran}{ran}

\newcommand{\mc}[1]{{\mathcal{#1}}}

\newcommand{\bb}[1]{{\mathbb{#1}}}

\begin{document}

\title{\textbf{Higher-order interlacing for matrix-valued meromorphic Herglotz functions}}
\author{Jakob Reiffenstein\thanks{Faculty of Mathematics, University of Vienna, Austria. The author was supported by the Austrian Science Fund [grant number I-4600].}}
\date{}

\maketitle
\vspace{-5ex}

\begin{abstract}
\noindent Scalar-valued meromorphic Herglotz-Nevanlinna functions are characterized by the interlacing property of their poles and zeros together with some growth properties. We give a characterization of matrix-valued Herglotz-Nevanlinna functions by means of a higher-order interlacing condition. As an application we deduce a matrix version of the classical Hermite-Biehler Theorem for entire functions.
\end{abstract}


\noindent © 2022. This manuscript version is made available under the CC-BY-NC-ND 4.0 license: https://creativecommons.org/licenses/by-nc-nd/4.0/ \\
The final publication is available at Elsevier via https://doi.org/10.1016/j.jmaa.2022.126260


\section{Introduction}
\noindent Interlacing patterns of the zeros of two polynomials have occurred already in the middle of the 19th century. The work of Sturm, Cauchy, and Hermite \cite{hermite} established a connection between the stability\footnote{A polynomial is called \textit{stable} if all of its zeros have negative real part.} of a polynomial $p$ in terms of interlacing of the real zeros of two associated polynomials, now known as the Hermite-Biehler Theorem. This criterion can be found in different versions in the work of Biehler \cite{biehler} and Hurwitz \cite{hurwitz}. Hurwitz in fact mentions the successful application of a related result during the construction of the turbine system of the Swiss bathing resort Davos. \\
In modern times, the Hermite-Biehler Theorem has been developed further and there are versions for polynomials with complex coefficients \cite{ho-da-ba-c}, multivariate polynomials \cite{wagner}, for entire functions \cite{levin}, and polynomials with $l$ zeros in the left half-plane and $r$ zeros in the right half-plane \cite{ho-da-ba}. The geometry of zeros of a polynomial is relevant in combinatorics \cite{braenden} as well as in spectral theory \cite{moeller-piv}, \cite{piv-wo} and in control theory \cite{si-da-ba}, where Kharitonov's Theorem \cite{kharitonov} is a prominent consequence of the Hermite-Biehler Theorem.  Moreover, the class of Hermite-Biehler functions is essential for de Branges' theory of Hilbert spaces of entire functions \cite{debranges}. Let us recall:

\begin{theo}[Hermite-Biehler]
\label{classichb}
Let $A$ and $B$ be nonzero polynomials with real coefficients. Set $p=A+iB$. Then all zeros of $p$ belong to the open lower half-plane $\bb{C}_-$, if and only if the following two conditions hold:
\begin{enumerate}
\item[$(i)$] The zeros of $A$ and $B$ are all real, simple, and interlace.
\item[$(ii)$] There is $x \in \mathbb{R}$ such that $A'(x)B(x)-B'(x)A(x)>0$.
\end{enumerate}
\end{theo}

\noindent A generalization of this theorem to entire functions is due to Levin, Meiman and Naimark and can be found in Chapter VII of Levin's book \cite{levin}. \\
\noindent Our motivation for the present work was to find a version of the Hermite-Biehler Theorem for matrix-valued functions. Items $(i)$, $(ii)$ in \Cref{classichb} essentially mean that $q:=\frac{A}{B}$ is a meromorphic Herglotz function, i.e., a function meromorphic on $\mathbb{C}$ taking real values on $\mathbb{R}$ that maps the open upper half-plane $\mathbb{C}_+$ to $\mathbb{C}_+ \cup \mathbb{R}$. The task of extending the theorem to matrix-valued functions thus boils down to investigating the pattern of zeros and poles of matrix-valued meromorphic Herglotz functions. For classical results on matrix-valued (not necessarily meromorphic) Herglotz functions, we refer the reader to \cite{gesztesy-tsekanovskii} for a general overview, and to \cite{alpay-tsekanovskii}, \cite{gkmt01} for the representation of matrix-valued and operator-valued Herglotz functions in terms of transfer matrices of certain linear stationary dynamical systems. \\

\noindent In the article at hand, after some preliminaries, we discuss in \Cref{interlacing_section} a generalized interlacing condition and its significance for scalar-valued meromorphic functions. In \Cref{sec_sufficiency}, we first prove that the zeros and poles of the determinant of any matrix-valued meromorphic Herglotz function are interlacing in this generalized sense. We then continue with the main result of this paper, \Cref{sufficiency}. In this theorem we show that matrix-valued meromorphic Herglotz functions can be fully characterized through the generalized interlacing of the zeros and poles of all principal minors.\\
\noindent In \Cref{hb_section}, we use this result to prove \Cref{final-matrix-hb}, which can be seen as a version of the Hermite-Biehler Theorem for matrix-valued functions.

\section{Preliminaries}
\noindent For a function $f: \Omega \to \bb{C}^{n \times n}$ defined on some set $\Omega \subseteq \bb{C}$, we set $f^{\#}(z):=f(\overline{z})^*$ whenever $\overline{z} \in \Omega$. We say that $f$ is $\#$-\textit{real} if $\Omega$ is symmetric w.r.t. the real axis, and $f=f^{\#}$ on $\Omega$. \\
The $*$-real and $*$-imaginary part of a (constant) complex matrix $M$ will be denoted by $\RE M:=\frac{1}{2}(M+M^*)$, while $\IM M:=\frac{1}{2i}(M-M^*)$. \\
\noindent The open lower and upper half-plane we will call $\bb{C}_{\pm}:=\{z \in \bb{C}|\pm\IM z >0\}$.
\begin{definition}
A Herglotz function is a function $Q: \bb{C} \setminus \bb{R} \to \mathbb{C}^{n \times n}$ that is holomorphic, \\ $\#$-real, and satisfies $\IM Q(z) \geq 0$ for all $z \in \bb{C}_+$. If $Q$ admits a meromorphic continuation to all of $\bb{C}$, we say that $Q$ is a \textit{meromorphic} Herglotz function. For convenience, the continuation is then denoted by $Q$ as well.
\end{definition}
\noindent Herglotz functions are often also called Nevanlinna functions, or Herglotz-Nevanlinna functions. Any piece of literature about them will contain the fundamental theorem giving a unique integral representation for every Herglotz function. As we will only need it for meromorphic Herglotz functions, the theorem will be given in a simpler version.

\begin{theo}
\label{matrixherg_int_rep}
Let $Q: \bb{C} \setminus \bb{R} \to \mathbb{C}^{n \times n}$. Then $Q$ is a meromorphic Herglotz function if and only if it admits, for all $z \in \bb C \setminus \bb R$, the following representation,
\begin{equation}
\label{matrix_discrete}
Q(z)=C+Dz+\sum_{j \in M} A_j \Big( \frac{1}{z_j-z}-\frac{z_j}{1+z_j^2} \Big)
\end{equation}
where $C=C^*$, $D \geq 0$, where $(A_j)_{j \in M}$ is a finite or infinite sequence of positive semi-definite
complex $n \times n$ matrices with
\begin{equation}
\label{A_j_summable}
\sum_{j \in M} \frac{\tr A_j}{1+z_j^2} < +\infty,
\end{equation}
and where $(z_j)_{j \in M}$ is a finite or infinite sequence of real numbers without finite limit points.
Moreover, 
\begin{enumerate}
\item[(i)] $C=\RE Q(i)$;
\item[(ii)] $D=\lim_{\eta \to +\infty} \frac{\IM Q(i\eta)}{\eta}=\lim_{\eta \to +\infty} \frac{Q(i\eta)}{i\eta}$;
\item[(iii)] Let $-\infty<c<d<+\infty$, where neither $c$ nor $d$ occur in the sequence $(z_j)_{z \in M}$. Then
\begin{equation*}
\sum_{\substack{j \in M \\ c<z_j<d}} A_j=\frac{1}{\pi} \lim_{\eta \searrow 0} \int_c^d \IM Q(x+i\eta) \, dx;
\end{equation*}
\item[(iv)] If all diagonal entries of $Q$ vanish, then $Q \equiv C$.
\item[(v)] If $\det Q \not\equiv 0$, then $\det Q(z) \neq 0$ for every nonreal $z$.
\end{enumerate}
\end{theo}
\vspace{7pt}
\noindent For a scalar-valued meromorphic Herglotz function $q \not\equiv 0$, there is more to say. Namely, the zeros and poles of $q$ interlace, which follows immediately from the fact that $q(x)$ is stricly increasing for real $x$. \\
\noindent Conversely, let $a_M,...,a_N$ and $b_{\hat{M}},...,b_{\hat{N}}$ with $M, \hat{M} \geq -\infty$ and $N, \hat{N} \leq +\infty$ be interlacing sequences, i.e., $b_k<a_k<b_{k+1}$, there is a unique, up to a multiplicative constant, meromorphic Herglotz function having zeros $\{a_k\}$ and poles $\{b_k\}$. To avoid technicalities, we only state the following result for the case $M=-\infty$ and $N=+\infty$, but a similar formula holds for every pair of interlacing sequences.

\begin{theo}[{\cite[Chapter VII, Theorem 1]{levin}}]
\label{herg_inf_poles}
Let the scalar-valued function $q$ be meromorphic on $\mathbb{C}$ and holomorphic on $\mathbb{C} \setminus \mathbb{R}$ such that the set of poles of $q$ is bounded neither from below nor from above.
Then $q$ is Herglotz if and only if
\begin{equation}
\label{herg_inf_poles_rep}
q(z)=c \frac{a_0-z}{b_0-z} \prod_{j \in \mathbb{Z} \setminus \{0\}} \frac{1-\frac{z}{a_j}}{1-\frac{z}{b_j}}, \quad z \in \mathbb{C} \setminus \mathbb{R},
\end{equation}
where $b_j<a_j<b_{j+1}$, $j \in \mathbb{Z}$, $a_{-1}<0<b_1$, and $c>0$. \\
For any such pair of interlacing sequences, the product converges uniformly on compact subsets of $\bb{C}$ not containing any of the points $b_k$.
\end{theo}

\section{The $n$-interlacing condition}
\label{interlacing_section}

\noindent In order to get an analog of \Cref{herg_inf_poles} for matrix-valued Herglotz functions, we need a higher-order version of the interlacing condition. We will formulate it in terms of divisor functions, where the function $\theta$ is called a divisor function if it is defined on a subset of $\bb{C}$ taking values in $\bb{Z}$ such that $\theta^{-1}(\bb{Z} \setminus \{0\})$ is a discrete set. For a scalar-valued meromorphic function $f \not\equiv 0$, denote by $\theta_f$ the divisor function of $f$, i.e., $\theta_f(z)=n$ if $z$ is a zero of multiplicity $n$ of $f$, $\theta_f(z)=-n$ if $z$ is a pole of multiplicity $n$ of $f$, and $\theta_f(z)=0$ otherwise.

\begin{definition}
\label{def_n-interlacing}
Let $\theta:\mathbb{R} \to \mathbb{Z}$ be a divisor function, and let $n \in \mathbb{N}$.
Then $\theta$ is called \emph{$\it{n}$-interlacing}, if for every $-\infty<a<b<+\infty$ we have
\begin{equation}
\label{n-interlacing}
\Big| \sum_{x \in (a,b)} \theta(x) \Big| \leq n.
\end{equation}
If $f$ is a meromorphic function on some open set $\Omega \supseteq \mathbb{R}$, we say that $f$ satisfies the \emph{$\it{n}$-interlacing condition} if $f \not\equiv 0$, all zeros and poles of $f$ are real and the function $\theta_f|_{\mathbb{R}}$ is $n$-interlacing.
\end{definition}

\noindent The following fact might also be of independent interest.

\begin{prop}
\label{colouring}
Let $\theta: \mathbb{R} \to \mathbb{Z}$ be a divisor function and let $n \in \mathbb{N}$. Then $\theta$ is $n$-interlacing if and only if there exist $1$-interlacing divisor functions $\theta_1,...,\theta_n$ such that $\theta=\sum_{j=1}^n \theta_j$.
\end{prop}
\begin{proof}
The backwards implication is evident. For the forward implication, assume that $\theta$ is $n$-interlacing.

\noindent Set
\begin{equation*}
\Theta(x)= \left\{ \begin{array}{ll}
\sum_{t \in [0,x)} \theta(t), &x \geq 0, \\
-\sum_{t \in (x,0)} \theta(t), &x < 0.
\end{array} \right.
\end{equation*}
Then $\Theta$ is a well defined step function because of $\theta$ having discrete support.
We use the notation $\Theta (x-)$ for $\lim_{t \nearrow x} \Theta(t)$ and $\Theta (x+)$ for $\lim_{t \searrow x} \Theta(t)$.
For $j \in \mathbb{Z}$, we define
\begin{equation*}
\theta_j(x):=\left\{ \begin{array}{ll}
1 & \text{if } \Theta(x+)> j \geq \Theta(x-), \\
-1 & \text{if } \Theta(x-) > j \geq \Theta(x+), \\
0 & \text{otherwise.}
\end{array} \right.
\end{equation*}
We show that each $\theta_j(x)$ is already $1$-interlacing. \\
Suppose that we are given $x<y$ such that $\theta_j(x)=\theta_j(y)=1$. Then $\Theta(x-) \leq j < \Theta(x+)$ and $\Theta(y-) \leq j < \Theta(y+)$. In particular, 
\begin{equation*}
\Theta(x+)>j \geq \Theta(y-).
\end{equation*}
Hence, $t_0:=\inf \{t > x:\Theta(t+) \leq j \}$ is well defined, and $t_0 \in (x,y)$ such that $\Theta(t_0-)>j \geq \Theta(t_0+)$. Therefore, $\theta_j(t_0)=-1$.
Analogously, we see that between points $\tilde{x}<\tilde{y}$ with $\theta_j(\tilde{x})=\theta_j(\tilde{y})=-1$, there is $\tilde{t_0}$ satisfying $\theta_j(\tilde{t_0})=1$. We conclude that $\theta_j$ is $1$-interlacing. \\
We check that $\theta=\sum_{j \in \bb{Z}} \theta_j$. Letting $x<y$, observe that
\begin{equation*}
|\Theta(y)-\Theta(x)| = \Big|\sum_{t \in (x,y)} \theta(t) \Big| \leq n.
\end{equation*}
Writing $j_-:= \min \{\Theta(x): x \in \mathbb{R}\}$ and $j_+ := \max \{\Theta(x): x \in \mathbb{R}\}$, it follows that $j_+-j_- \leq n$. However, from our definition of $\theta_j$ we can see that $\theta_j \equiv 0$ if $j \geq j_+$ or $j<j_-$, which means that all but at most $n$ of the functions $\theta_j$ vanish. \\
Now, the definition of $\Theta$ yields 
\begin{align*}
\theta(x) &=\Theta(x+)-\Theta(x-) \\
&=\left\{ \begin{array}{ll}
\big| \{j \in \mathbb{Z}: \Theta(x+)>j \geq \Theta(x-) \} \big|, &\Theta(x+)>\Theta(x-), \\
-\big| \{j \in \mathbb{Z}: \Theta(x-)>j \geq \Theta(x+) \} \big|, &\Theta(x+)<\Theta(x-), \\
0, &\Theta(x+)=\Theta(x-)
\end{array} \right. \\
&=\sum_{j \in \mathbb{Z}} \theta_j(x), \quad x \in \mathbb{R}.
\end{align*}
We already know that the number of non-vanishing $\theta_j$ is at most $n$. The proof is complete.
\end{proof}
  
\noindent Let $f$ be a scalar-valued meromorphic function on $\mathbb{C}_+$. We say that $f$ is of \textit{bounded type} if it can be written as the ratio $g/h$ of two functions $g,h$ defined on $\mathbb{C}_+$ that are analytic and bounded\footnote{This definition is equivalent to the one given in \cite[Section 3.1]{arovdym}. It differs from Definition 3.15 given in \cite{roro}. However, for holomorphic functions on a simply connected domain (such as $\mathbb{C}_+$), the two definitions are equivalent by \cite[Theorem 3.20]{roro}. In this case, $h$ may be chosen to not have any zeros.}. \\
For a function $f \not\equiv 0$ of bounded type, the \textit{mean type}
\begin{equation}
\tau_f:= \limsup_{\eta \to +\infty} \frac{\log |f(i\eta)|}{\eta}
\end{equation}
is finite.  \\
Notably, the restriction to $\mathbb{C}_+$ of a Herglotz function is always of bounded type, as is seen by composing with a M\"{o}bius transformation mapping the upper half-plane to the unit circle.

\begin{theo}
\label{colour_func}
Let $n \in \mathbb{N}$ and let $f \not\equiv 0$ be meromorphic on $\mathbb{C}$. Then $f$ satisfies the $n$-interlacing condition if and only if there exist scalar-valued meromorphic Herglotz functions $q_j \not\equiv 0$, $j \in \{1,...,n\}$, a complex constant $C$ with $|C|=1$, and an entire function $g$, such that $g(0)=0$ and
\begin{equation}
\label{func_n_interl}
f=Ce^g\prod_{j=1}^n q_j.
\end{equation}
If, additionally, $f$ is $\#$-real and $f|_{\mathbb{C}_+}$ is of bounded type, then $g=0$ and $C=\pm 1$.
\end{theo}
\begin{proof}
Suppose that $f$ can be represented as in (\ref{func_n_interl}), for an entire function $g$ and meromorphic Herglotz functions $q_j \not\equiv 0$, $j=1,...,n$. Then
\[
\theta_f=\theta_{e^g}+\sum_{j=1}^n \theta_{q_j}=\sum_{j=1}^n \theta_{q_j}.
\]
Using \Cref{herg_inf_poles}, we find that $\theta_{q_j}|_{\mathbb{R}}$ is $1$-interlacing for each $j \in \{1,...,n\}$. Hence, by \Cref{colouring}, $\theta_f|_{\mathbb{R}}$ is $n$-interlacing, i.e., $f$ satisfies the $n$-interlacing condition. \\
On the other hand, if $f$ satisfies the $n$-interlacing condition, \Cref{colouring} provides $1$-interlacing functions $\theta_1,...,\theta_n$ such that $\theta_f|_{\mathbb{R}}=\sum_{j=1}^n \theta_n$. Let also $\theta_j^{-1}(\{1\})=\{a_{j,M_j},...,a_{j,N_j} \}$ and $\theta_j^{-1}(\{-1\})=\{b_{j,\hat{M_j}},...,b_{j,\hat{N_j}}\}$, where $M_j,\hat{M_j} \geq -\infty$, $N_j,\hat{N_j} \leq +\infty$ and
\begin{align*}
&\cdots < b_{j,k} < a_{j,k} < b_{j,k+1} < a_{j,k+1} < \cdots
\end{align*}
We want to define a function having its zeros at the points $a_{j,k}$ and its poles at the points $b_{j,k}$, in the same way as in (\ref{herg_inf_poles_rep}). Providing full generality requires some technicalities, however. Let
\begin{equation*}
r_j=\prod_{k \in \bb{Z}} m_{j,k}
\end{equation*}
where $m_{j,k}$ is the meromorphic function on $\mathbb{C}$ defined by
\begin{align*}
m_{j,k}(z):=\left\{ \begin{array}{ll}
\frac{1-z/a_{j,k}}{1-z/b_{j,k}}, & \hat{M_j} \leq k \leq N_j, \quad a_{j,k} \neq 0, b_{j,k} \neq 0;  \\
\frac{a_{j,k}-z}{b_{j,k}-z}, & \hat{M_j} \leq k \leq N_j, \quad ( a_{j,k} = 0 \vee b_{j,k} = 0);  \\
a_{j,k}-z, & k=M_j<\hat{M_j}; \\
\frac{1}{b_{j,k}-z}, & k=\hat{N_j}>N_j; \\
1, & \text{else.}
\end{array} \right.
\end{align*}
By \Cref{herg_inf_poles}, the product $r_j$ converges locally uniformly, and either $r_j$ or $-r_j$ is a Herglotz function. Letting $q_j:=\pm r_j$ such that $q_j$ is Herglotz, the function $h:=\frac{f}{q_1 \cdots q_n}$ is entire and does not have any zeros, and therefore can be written as $Ce^g$, with $C \neq 0$ and some entire function $g$ satisfying $g(0)=0$. Since $|C|q_1$ is Herglotz as well, one may choose $C$ such that $|C|=1$. Hence, representation (\ref{func_n_interl}) holds.  \\
\noindent Suppose now that $f$ is also $\#$-real and that $f|_{\mathbb{C}_+}$ is of bounded type. Since, for every $j \in \{1,...,n\}$, the function $q_j$ is $\#$-real and of bounded type, the same goes for $h=Ce^g$. In view of \cite[Problem 34, Lemma 2]{debranges}, we have $h(z)=C \exp (-az^2-ibz)$ for some constants $a,b$ and all $z \in \mathbb{C}$. Using implication $(ii) \Rightarrow (i)$ of \cite[Theorem 6.17]{roro}, we obtain that $h$ is of exponential type\footnote{An entire function $F$ is \textit{of exponential type} if $\sigma_F:=\limsup_{|z| \to +\infty} \frac{\log |F(z)|}{|z|}<\infty$.} and that $\int_{\bb R} \frac{\max \{0, \log |h(t)|\}}{1+t^2} \, dt <+\infty$. It is now easy to see that $a=b=0$. As $h$ is $\#$-real, it follows that $C$ is real and thus equal to $\pm 1$.
\end{proof}

\noindent The relevance of the next statement is the characterization of the case that the functions $q_1,...,q_n$ can be chosen to be meromorphic.

\begin{theo}
\label{charprodherg}
Let $f$ be meromorphic on $\mathbb{C}$, nonconstant, and $\#$-real. For $n \in \mathbb{N}$, the following conditions are equivalent:
\begin{enumerate}
\item[$(i)$] There exist meromorphic Herglotz functions $q_1,...,q_n$ such that $f=\prod_{j=1}^n q_j$;
\item[$(ii)$] There exists a holomorphic logarithm $u$ of $f|_{\mathbb{C}_+}$ with $\IM u(z) \in (0,n\pi)$ for every $z \in \mathbb{C}_+$.
\end{enumerate}
\end{theo}
\begin{proof}
Suppose $(i)$. If $f=\prod_{j=1}^n q_j$, then $u:=\sum_{j=1}^n \log q_j|_{\mathbb{C}_+}$ is a logarithm of $f$ with $\IM u(z) \in (0,n\pi)$, $z \in \mathbb{C}_+$. \\
\noindent We show that $(i)$ follows from $(ii)$. Consider any finite interval $(a,b)$ and assume that neither $a$ nor $b$ is a pole or a zero of $f$. Set $x_0:=\frac{a+b}{2}$ and $r:=\frac{a-b}{2}$. Define $\gamma_+(t):=x_0+re^{it}$ for $t \in [0,\pi]$ and $\gamma_-(t):=x_0+re^{it}$ for $t \in [\pi,2\pi]$.
Then 
\begin{equation}
\label{logint}
\sum_{x \in (a,b)} \theta_f(x)=\frac{1}{2\pi i} \bigg( \int_{\gamma_+} \frac{f'(z)}{f(z)} \, dz + \int_{\gamma_-} \frac{f'(z)}{f(z)} \, dz \bigg).
\end{equation}
Observe that $f=f^{\#}$ implies that $u^{\#}$ is a holomorphic logarithm of $f^{\#}|_{\mathbb{C}_-}$ satisfying $\IM g^{\#}(z) \in (-n\pi,0)$, $z \in \mathbb{C}_-$. Thus, the modulus of the imaginary parts of both integrals in (\ref{logint}) is bounded by $n\pi$ each. So, $|\sum_{x \in (a,b)} \theta_f(x)| \leq n$. Applying \Cref{colour_func} yields meromorphic Herglotz functions $q_1,...,q_n$, a real constant $C$, and an entire $\#$-real function $g$ such that $f=Ce^g\prod_{j=1}^n q_j$. Again, each $q_j$ has a holomorphic logarithm in the upper half-plane, which leads to $g|_{\mathbb{C}_+}=\log C+u-\sum_{j=1}^n \log q_j$. Therefore, $|\IM g|$ is bounded in $\mathbb{C}_+$. Since $g$ is $\#$-real, we get that $|\IM g|$ is also bounded in $\mathbb{C}_-$ and thus in all of $\mathbb{C}$. By the Liouville Theorem, $g$ is constant, so $f=\tilde{C} \prod_{j=1}^n q_j$ for some real constant $\tilde{C}$. If $\tilde{C} \geq 0$, then $f=(\tilde{C}q_1)\prod_{j=2}^n q_j$, which is the desired representation.
This leaves the case $\tilde{C}<0$. First, we write $f=-(-\tilde{C}q_1)\prod_{j=2}^n q_j$, where $-\tilde{C}q_1$ is Herglotz. Then the function $\tilde{u}:=i\pi+\log (-\tilde{C}q_1)+\sum_{j=2}^n \log q_j$ is another holomorphic logarithm of $f$ and satisfies $\IM \tilde{u}(z) \in (\pi,(n+1)\pi)$ for $z$ in $\mathbb{C}_+$. Observing that $u=2r\pi i+\tilde{u}$ for some $r \in \mathbb{Z}$, we get that $(2r+1)\pi < \IM u(z) < (2r+1+n)\pi$ for all $z \in \mathbb{C}_+$. At the same time, we know that $0 < \IM u < n\pi$. In total, we either have $\IM u \in (0,(n-1)\pi)$ (if $r<0$) or $\IM u \in (\pi,n\pi)$. In both cases, from what was proven so far, we conclude that there exist meromorphic Herglotz functions $r_1,...,r_{n-1}$ such that $f=\prod_{j=1}^{n-1} r_j$ or $f=-\prod_{j=1}^{n-1} r_j$. Setting $\delta \equiv 1$ and $\delta \equiv -1$, respectively, leads to $f=\delta\prod_{j=1}^{n-1} r_j$, which proves the theorem because $\delta$ is a meromorphic Herglotz function.
\end{proof}

\section{Characterization of matrix-valued meromorphic Herglotz functions}
\label{sec_sufficiency}

\noindent Our goal in this section is to prove a statement of the form ``$Q$ is an $n \times n$-matrix-valued meromorphic Herglotz function if and only if its zeros and poles are $n$-interlacing''. It is not yet clear how to understand this, but we see in the next lemma that we should look at the determinant of $Q$.

\begin{lemma}
\label{prodofscalarherg}
Let $Q$ be an $n \times n$-matrix-valued Herglotz function. Then there exist scalar-valued Herglotz functions $q_1,...,q_n$ such that $\det Q=q_1\cdots q_n$. If $Q$ is meromorphic on $\mathbb{C}$, then $q_1,...,q_n$ can be chosen to be meromorphic as well.
\end{lemma}
\begin{proof}
The proof proceeds by induction on $n$. Since the assertion is evident for $n=1$,
only the induction step is to be done. Suppose that the assertion holds for every $n \times n$-matrix-valued Herglotz function, and let $Q$ be a $(n+1) \times (n+1)$-matrix-valued Herglotz function. If $Q \equiv C$ (with self-adjoint $C$), then there is nothing to be proven. The same holds for the case that  $\det Q \equiv 0$. Otherwise, by item $(v)$ in \Cref{matrixherg_int_rep}, $\det Q(z) \neq 0$ for all $z \in \mathbb{C}_+$. Letting $Q_{(j)}$, $j=1,...,n+1$, be the submatrix of $Q$ obtained by deleting the $j$-th row and column from $Q$, we can write
\begin{align*}
-Q^{-1}=\left( \begin{array}{cccc}
-\frac{\det Q_{(1)}}{\det Q} & * & * & * \\
*& -\frac{\det Q_{(2)}}{\det Q} & *& *\\
*&*& \ddots &* \\
*&*&*& -\frac{\det Q_{(n+1)}}{\det Q}
\end{array} \right)
\end{align*} 
Because $-Q^{-1}$ is Herglotz, all of its diagonal entries $-\frac{\det Q_{(j)}}{\det Q}$ are scalar-valued Herglotz functions. Because we assumed $Q$ not to be constant, and because of \Cref{matrixherg_int_rep}, $(iv)$, there exists an integer $j_0$ such that $1 \leq j_0 \leq n+1$ and $q:=-\frac{\det Q_{(j_0)}}{\det Q}$ is not identically zero. Note that $-q^{-1}$ is also a scalar Herglotz function, which implies the assertion by the fact that $Q_{(j_0)}$ is a matrix-valued Herglotz function and by
\begin{equation}
\det Q=\det Q_{(j_0)}[-q^{-1}]=q_1\cdots q_n \cdot [-q^{-1}].
\end{equation}
\end{proof}

\noindent We continue with a reminder on linear algebra. Fix $n \in \mathbb{N}$ and let $1 \leq m \leq n$. Set
\[
I_{n,m}:=\big\{(i_1,...,i_m) \in \mathbb{N}^m \mid 1 \leq i_1<...<i_m \leq n \big\}.
\]
For $\mathbf{i} \in \bigcup_{m=1}^n I_{n,m}$, let $|\mathbf{i}|$ be the length of $\mathbf{i}$, i.e., $|\mathbf{i}|$ is such that $\mathbf{i} \in I_{n,|\mathbf{i}|}$.

\noindent Let $M=(m_{ij}) \in \mathbb{C}^{n \times n}$, $m \in \{1,...,n\}$, and $\mathbf{i} \in I_{n,m}$. We define the matrix $M_{\mathbf{i}}$ by
\begin{align*}
M_{\mathbf{i}}:=\left( \begin{array}{cccc}
m_{i_1i_1} & m_{i_1i_2} & \cdots & m_{i_1i_m} \\
m_{i_2i_1} & m_{i_2i_2} & \cdots & m_{i_2i_m} \\
\vdots & \vdots & \ddots & \vdots \\
m_{i_mi_1} & m_{i_mi_2} & \cdots & m_{i_mi_m} 
\end{array} \right) \, \in \mathbb{C}^{m \times m}.
\end{align*}
A matrix $X$ is called a \emph{principal submatrix} of $M$ if there exists $1 \leq m \leq n$ and $\mathbf{i} \in I_{n,m}$, such that $X=M_{\mathbf{i}}$. A \emph{principal minor} of $M$ is the determinant of a principal submatrix of $M$. \\

\noindent In the following lemma, the restriction to $\mathbb{C}_+$ of $Q$ is assumed to be a matrix-valued function of bounded type.  This means that all entries of $Q|_{\mathbb{C}_+}$ are scalar-valued functions of bounded type.

\begin{lemma}
\label{sum_rep}
Let $Q$ be an $n \times n$-matrix-valued function that is meromorphic on $\mathbb{C}$ and such that all poles of $Q$ are real. Suppose that $Q$ is $\#$-real, $Q|_{\mathbb{C}_+}$ is of bounded type, and $\limsup_{\eta \to +\infty} \frac{\|Q(i\eta)\|}{\eta} <+\infty$. Let $\det Q_{\mathbf{i}}$ be $|\mathbf{i}|$-interlacing for all $\mathbf{i} \in I_{n,1} \cup I_{n,2}$. For every $\mathbf{i} \in I_{n,2}$, if $\det Q_{\mathbf{i}}$ has a pole of multiplicity $2$, assume that there is $z^* \in \bb R$ such that $\lim_{z \to z^*} (z^*-z)^2 \det Q_{\mathbf{i}}(z) >0$. Then $Q$ can be represented as
\begin{equation}
\label{hergrep_proof}
Q(z)=C+Dz+\sum_{j \in M} A_j \bigg(\frac{1}{z_j-z}-\frac{z_j}{1+z_j^2} \bigg), \quad z \in \bb C \setminus \bb R,
\end{equation}
where $\{z_j \mid j \in M \} \subseteq \mathbb{R}$ is the set of poles of $Q$, and $C,D,A_j$, $j \in M$, are Hermitian matrices and where
\begin{align}
\label{weighted_summability}
\sum_{j \in M} \frac{|A_{j;kl}|}{1+z_j^2} < + \infty \hspace{10ex} \textrm{for all } k,l \in \{1,...,n\}
\end{align}holds true.
\end{lemma}
\begin{proof}
We start by observing that, if (\ref{hergrep_proof}) holds, then $C,D,A_j$, $j \in M$, are all Hermitian. Clearly, $C=\RE Q(i)$ and $D=\lim_{\eta \to +\infty} \frac{\IM Q(i\eta)}{\eta}$ are Hermitian. Moreover, for every $j \in M$,
\[
A_j^*=\big[\lim_{z \to z_j} (z_j-z) Q(z) \big]^*=\lim_{z \to z_j} (z_j-\overline{z})Q(\overline{z})=A_j.
\]
We will show (\ref{hergrep_proof}) entrywise, i.e., for any $k, l \in \{1,...,n\}$, 
\begin{equation*}
Q_{kl}(z)=C_{kl}+D_{kl}z+\sum_{j \in M} A_{j;kl} \bigg(\frac{1}{z_j-z}-\frac{z_j}{1+z_j^2} \bigg), \quad z \in \mathbb{C} \setminus \mathbb{R}.
\end{equation*}
If $k=l$, this representation is possible: By assumption, $Q_{kk}$ is $1$-interlacing, and \Cref{colour_func} shows that $Q_{kk}$ or $-Q_{kk}$ is a scalar-valued meromorphic Herglotz function. From \Cref{matrixherg_int_rep} we can thus infer
\begin{align}
\sum_{j \in M} \frac{|A_{j;kk}|}{1+z_j^2} < +\infty \hspace{10ex} \text{for all } k \in \{1,...,n\}.
\end{align}
Assume that $k<l$ and set $\mathbf{i}=(k,l) \in I_{n,2}$. We first show that all poles of $Q_{kl}$ are simple. Our assumption together with \Cref{colour_func} yields Herglotz functions $r_1,r_2$ such that $\det Q_{\mathbf{i}}=\pm r_1r_2$. Because of $Q_{lk}=Q_{kl}^{\#}$, we have $\theta_{Q_{kl}}=\theta_{Q_{lk}}$ and thus, for every $x \in \mathbb{R}$,
\begin{align*}
\big|\theta_{Q_{kl}}(x)\big| &=\frac{1}{2} \big|\theta_{Q_{kl}Q_{lk}}(x) \big|=\frac{1}{2} \big|\theta_{Q_{kk}Q_{ll} \mp r_1r_2}(x) \big| \leq \\
&\leq \frac{1}{2}\max \big\{|\theta_{Q_{kk}}(x)+\theta_{Q_{ll}}(x)|,|\theta_{r_1}(x)+\theta_{r_2}(x)| \big\} \leq 1.
\end{align*}
Let $j \in M$, and let $A_j$ be the residue of $Q$ at $z_j$. Then 
\begin{align*}
\det (A_j)_{\mathbf{i}}&=\det \big[\lim_{z \to z_j} (z_j-z)Q_{\mathbf{i}}(z) \big] \\
&=\lim_{z \to z_j} (z_j-z)^2 \det Q_{\mathbf{i}}(z)=\pm \underbrace{\big[\lim_{z \to z_j} (z_j-z)r_1(z) \big]\big[\lim_{z \to z_j} (z_j-z)r_2(z) \big]}_{\geq 0}.
\end{align*}
Suppose that $\det (A_j)_{\mathbf{i}} <0$. Then $\det Q_{\mathbf{i}}=-r_1r_2$, and hence $\det (A_j)_{\mathbf{i}} \leq 0$ for all $j \in M$. On the other hand, $\det Q_{\mathbf{i}}$ has a pole of multiplicity $2$ at $z_j$, and we find $z^* \in \mathbb{R}$ such that $\lim_{z \to z^*} (z^*-z)^2 \det Q_{\mathbf{i}}(z) >0$. This is a contradiction. Hence, for all $j \in M$, we have $\det (A_j)_{\mathbf{i}} \geq 0$ and thus
\begin{equation*}
|A_{j;kl}|=\big|A_{j;kl}A_{j;lk}\big|^{\frac{1}{2}} \leq \big|A_{j;kk}A_{j;ll}\big|^{\frac{1}{2}}.
\end{equation*}
Now \Cref{matrixherg_int_rep} and (\ref{A_j_summable}) imply that
\begin{align*}
\sum_{j \in M} \frac{|A_{j;kl}|}{1+z_j^2} \leq \sum_{j \in M} \bigg(\frac{|A_{j;kk}|}{1+z_j^2}\bigg)^{\frac{1}{2}} \bigg(\frac{|A_{j;ll}|}{1+z_j^2}\bigg)^{\frac{1}{2}} \leq \bigg(\sum_{j \in M} \frac{|A_{j;kk}|}{1+z_j^2} \bigg)^{\frac{1}{2}} \bigg(\sum_{j \in M} \frac{|A_{j;ll}|}{1+z_j^2} \bigg)^{\frac{1}{2}} < +\infty.
\end{align*}
This ensures uniform convergence of 
\begin{equation*}
\psi_{kl}(z):=\sum_{j \in M} A_{j;kl} \bigg(\frac{1}{z_j-z}-\frac{z_j}{1+z_j^2} \bigg)=\sum_{j \in M} A_{j;kl} \frac{1+z_jz}{(z_j-z)(1+z_j^2)}
\end{equation*}
for $z$ in any compact subset of $\mathbb{C}$ not containing any of the points $z_j$. If we write $A_{j;kl}=A_{j;kl}^{r,+}-A_{j;kl}^{r,-}+i\Big(A_{j;kl}^{i,+}-A_{j;kl}^{i,-}\Big)$, such that the four numbers on the right side of this equality are nonnegative, we see that 
\begin{equation*}
\psi_{kl}=\psi_{kl}^{r,+}-\psi_{kl}^{r,-}+i\Big(\psi_{kl}^{i,+}-\psi_{kl}^{i,-}\Big)
\end{equation*}
is a linear combination of Herglotz functions, and thus is a function of bounded type. \\
\noindent The function $g_{kl}:=Q_{kl}-\psi_{kl}$ can be continued to an entire function, denoted by $g_{kl}$ as well. We would like to show that it is of exponential type $\sigma_{g_{kl}}=0$. Let $g_{kl}^+$ be the restriction of $g_{kl}$ to $\mathbb{C}_+$. 
From the assumption on $Q$ we obtain that $g_{kl}^+$ is of bounded type and satisfies $\limsup_{\eta \to +\infty} \frac{|g_{kl}^+(i\eta)|}{\eta} <+\infty$. In addition, $\tau_{g_{kl}^+} \leq 0$, as otherwise we would find $\delta>0$ and a sequence $\eta_n \nearrow +\infty$ such that $\log |g_{kl}^+(i\eta_n)| \geq \delta \eta_n$ for all $n \in \mathbb{N}$, leading to the contradiction 
\[
\limsup_{\eta \to +\infty} \frac{|g_{kl}^+(i\eta)|}{\eta} \geq \limsup_{n \to \infty} \frac{|g_{kl}^+(i\eta_n)|}{\eta_n} \geq \lim_{n \to \infty} \frac{e^{\delta \eta_n}}{\eta_n} = +\infty.
\]
Note that $g_{kl}^-:=g_{kl}^{\#}|_{\mathbb{C}_+}=g_{lk}^+$ is also of bounded type and such that $\tau_{g_{kl}^-} \leq 0$. Using implication $(ii) \Rightarrow (i)$ of \cite[Theorem 6.17]{roro}, we obtain that $g_{kl}$ is of exponential type. Moreover, by \cite[Theorem 6.18]{roro}, $\tau_{g_{kl}^+}+\tau_{g_{kl}^-} \geq 0$, showing that both mean types are equal to $0$. Now the second assertion of \cite[Theorem 6.18]{roro} yields
\[
\sigma_{g_{kl}}=\max \big\{\tau_{g_{kl}^+},\tau_{g_{kl}^-} \big\}=0.
\]
Hence, the entire function $\tilde{g}_{kl}$ defined by $\tilde{g}_{kl}(z):=\frac{g_{kl}(z)-g_{kl}(0)}{z}$ for $z \in \mathbb{C}$ is of exponential type $0$ as well. In addition to that, it is bounded on the imaginary axis, which means that $\tilde{g}_{kl}$ is constant by the Phragm\'en-Lindel\"{o}f principle (\cite{levin}). This implies $g_{kl}(z)=C_{kl}+D_{kl}z$ for all $z \in \mathbb{C}$ with some constants $C_{kl},D_{kl}$, which proves representation (\ref{hergrep_proof}). 
\end{proof}

\noindent We are now ready to state our main theorem. In essence, it states that a $\#$-real matrix-valued meromorphic function $Q$ is Herglotz if and only if its principal minors all satisfy the interlacing properties of their respective sizes.\\
\noindent However, it is necessary to impose some additional assumptions on the growth of the function $Q$. 

\begin{theo}
\label{sufficiency}
Let $Q$ be an $n \times n$-matrix-valued function that is meromorphic on $\mathbb{C}$ and such that all poles of $Q$ are real. Suppose that $Q$ is $\#$-real, $Q|_{\mathbb{C}_+}$ is of bounded type, and $\limsup_{\eta \to +\infty} \frac{\|Q(i\eta)\|}{\eta} <+\infty$. Then the following statements are equivalent:
\begin{enumerate}
\item[$(i)$] $Q$ is Herglotz;
\item[$(ii)$] For every $m \in \{1,...,n\}$ and $\mathbf{i} \in I_{n,m}$, the function $f_{\mathbf{i}}:=\det Q_{\mathbf{i}}$ satisfies at least one of the following properties: 
\begin{enumerate}
\item[a.] $f_{\mathbf{i}}$ satisfies the $m$-interlacing condition (cf. \Cref{def_n-interlacing}). If $f_{\mathbf{i}}$ has a pole of multiplicity $m$, there is $z^* \in \mathbb{R}$ such that $\lim_{z \to z^*} [(z^*-z)^mf_{\mathbf{i}}(z)]>0$. The limit $\lim_{\eta \to +\infty} \frac{f_{\mathbf{i}}(i\eta)}{(i\eta)^m}$, which exists by the proof of \Cref{colour_func}, is nonnegative.
\item[b.] $f_{\mathbf{i}}$ can be represented as the product of $m$ scalar-valued meromorphic Herglotz functions;
\item[c.] $f_{\mathbf{i}}|_{\mathbb{C}_+}$ has a holomorphic logarithm $u$ satisfying $\IM u(z) \in (0,m\pi)$ for all $z \in \mathbb{C}_+$;
\end{enumerate}
\item[$(iii)$] $\det Q_{\mathbf{i}}$ is $|\mathbf{i}|$-interlacing for all $\mathbf{i} \in I_{n,1} \cup I_{n,2}$. For every $\mathbf{i} \in I_{n,2}$, if $\det Q_{\mathbf{i}}$ has a pole of multiplicity $2$, there is $z^* \in \bb R$ such that $\lim_{z \to z^*} [(z^*-z)^2 \det Q_{\mathbf{i}}(z)] >0$. All residues of $Q$ are negative semi-definite, and $D:=\lim_{\eta \to +\infty} \big(\frac{\IM Q(i\eta)}{\eta}\big)$, which exists by \Cref{sum_rep}, is positive semi-definite.
\end{enumerate}
If $(i)-(iii)$ hold, then in $(ii)$, all of $a.-c.$ hold whenever $f_{\mathbf{i}}$ is not constant.
\end{theo}

\begin{proof}
\noindent $(i) \Rightarrow (ii)$: \\
Any principal submatrix $Q_{\mathbf{i}}$ of $Q$ is again a Herglotz function. Hence \textit{b.} holds because of \Cref{prodofscalarherg}. By \Cref{colour_func}, $a.$ is fulfilled, too. If $\mathbf{i} \in \bigcup_{m=1}^n I_{n,m}$ is such that $f_{\mathbf{i}}$ is not constant, then $c.$ also holds. \\[0.7ex]

\noindent $(i) \Rightarrow (iii)$: \\
In view of \Cref{matrixherg_int_rep}, the residues of $Q$ are all negative semi-definite, and $D$ is positive semi-definite. The remaining properties follow as in the proof of the implication $(i) \Rightarrow (ii)$. \\[0.7ex]

\noindent $(iii) \Rightarrow (i)$: \\
\Cref{sum_rep} shows that $Q$ has the representation (\ref{hergrep_proof}), where (\ref{weighted_summability}) is fulfilled and where $D \geq 0$ by our assumption. Since the residue of $Q$ at $z_j$ is $-A_j$, we also get that $A_j \geq 0$. Hence, (\ref{A_j_summable}) follows. By \Cref{matrixherg_int_rep}, $Q$ is Herglotz. \\[0.7ex]

\noindent $(ii) \Rightarrow (i)$: \\
Representation (\ref{hergrep_proof}) and (\ref{weighted_summability}) again follow from \Cref{sum_rep}. We fix $m \in \{1,...,n\}$ and $\mathbf{i} \in I_{n,m}$. \Cref{colour_func} yields Herglotz functions $r_1,...,r_m$ such that $f_{\mathbf{i}}=\det Q_{\mathbf{i}}=\pm \prod_{k=1}^m r_k$. For any $j \in M$, we have
\begin{align}
\label{res_reform}
\det (A_j)_{\mathbf{i}}=\det \big[\lim_{z \to z_j} (z_j-z)Q_{\mathbf{i}}(z) \big]
=\lim_{z \to z_j} [(z_j-z)^m f_{\mathbf{i}}(z)]=\pm \lim_{z \to z_j} \big[(z_j-z)^m \prod_{k=1}^m r_k(z)\big]
\end{align}
and
\begin{equation}
\label{mpolpos}
\lim_{z \to z_j} \big[(z_j-z)^m \prod_{k=1}^m r_k(z)\big]= \prod_{k=1}^m \big[\lim_{z \to z_j}(z_j-z)r_k(z) \big] \geq 0.
\end{equation}
If $b.$ is satisfied, this leads to $\det (A_j)_{\mathbf{i}}\geq 0$. The same is true if $c.$ holds, as $c.$ implies $b.$ Let us consider the case that $a.$ holds, and there is $j \in M$ such that $\det (A_j)_{\mathbf{i}} < 0$. By (\ref{res_reform}) and (\ref{mpolpos}), $f_{\mathbf{i}}=- \prod_{k=1}^m r_k$ and hence even $\det (A_j)_{\mathbf{i}} \leq 0$ for all $j \in M$. However, as $f_{\mathbf{i}}$ has a pole of multiplicity $m$ at $z_j$, we can find $z^* \in \mathbb{R}$ with $\lim_{z \to z^*} [(z^*-z)^m f_{\mathbf{i}}(z)] >0$, which is a contradiction. This shows that $\det (A_j)_{\mathbf{i}} \geq 0$ has to hold for every $m \in \{1,...,n\}$, $\mathbf{i} \in I_{n,m}$, and every $j \in M$. \\
Applying Sylvester's criterion (see \cite{meyer}, equations (7.6.9)-(7.6.12)), we obtain $A_j \geq 0$ for all $j \in M$. \\
Finally, it needs to be shown that $D \geq 0$. As a consequence of the proof of \Cref{sum_rep}, the matrix $\lim_{\eta \to +\infty} \frac{Q(i\eta)}{i\eta}$ exists and coincides with $D$. Choose $\mathbf{i}=(i_1,...,i_m)$ as above. Then $f_{\mathbf{i}}=\det Q_{\mathbf{i}}$ satisfies one of $a.-c$. If it satisfies $a.$, then by assumption
\begin{equation*}
\det D_{\mathbf{i}}=\lim_{\eta \to +\infty} \frac{f_{\mathbf{i}}(i\eta)}{(i\eta)^m} \geq 0.
\end{equation*}
Otherwise, $b.$ is satisfied. Therefore, $\tilde{f_{\mathbf{i}}}(z):=f_{\mathbf{i}}(-\frac{1}{z})$ still satisfies $b.$, which leads to \\ $\det D_{\mathbf{i}}=\lim_{\eta \to +\infty} \frac{f_{\mathbf{i}}(i\eta)}{(i\eta)^m}=\lim_{z \to 0} \Big((-z)^m \tilde{f_{\mathbf{i}}}(z)\Big) \geq 0$. By Sylvester's criterion, $D \geq 0$. The application of \Cref{matrixherg_int_rep} completes the proof.
\end{proof}

\section{The Hermite-Biehler Theorem for de Branges matrices}
\label{hb_section}
\noindent When extending the classical Hermite-Biehler Theorem from polynomials to entire functions, in order for the zeros and poles of $A$ and $B$ to interlace, it is not anymore sufficient that the zeros of $E$ lie in the lower half-plane. Instead, one introduces the class HB as the set of all entire functions that have no real zeros, and satisfy
\begin{equation}
|E(z)|>|E^{\#}(z)|, \quad z \in \bb{C}_+,
\end{equation}
see e.g. (\cite{levin}).
The most straightforward way to define a matrix-valued Hermite-Biehler class HB$_n$ is to require its elements $E:\bb{C} \to \bb{C}^{n \times n}$ to be entire and satisfy $\det E \not\equiv 0$ and
\begin{equation}
\label{hb_ineq}
E(z)E(z)^*-E^{\#}(z)E^{\#}(z)^* \geq 0, \quad z \in \bb{C}_+.
\end{equation}
Sometimes it is more convenient to set $s_E:=E^{-1}E^{\#}|_{\mathbb{C}_+}$ and require it to belong to the Schur class
\begin{equation}
\mc{S}^{n \times n} := \Big\{s: \bb{C}_+ \to \bb{C}^{n \times n} \text{ holomorphic}\,\mid\, I-s(z)s(z)^* \geq 0,\, z \in \bb{C}_+\Big\}.
\end{equation}
In this case, for real $x$ we get that $E(x)E(x)^*-E(x)^*E(x)$ is positive semi-definite while having zero trace, which is only possible for the zero matrix. Hence, $s_E$ even belongs to the subclass of \textit{inner} Schur functions
\begin{equation}
\mc{S}_{in}^{n \times n} := \Big\{s \in \mc{S}^{n \times n} \,\mid\,I-s(x)s(x)^*=0 \, \text{ for a.e. } x \in \bb{R} \Big\}.
\end{equation}
This set is well defined since it can be shown that every $f \in \mc{S}^{n \times n}$ has nontangential boundary values at almost all points $x \in \bb{R}$.
\medskip

\noindent We introduce a generalization of the class HB$_n$, adopting the terminology from \cite{arovdym}, Section 5.10.

\begin{definition}
A meromorphic function 
\begin{equation}
\mc{E}=[E_-,E_+]: \bb{C}_+ \to \bb{C}^{n \times 2n}
\end{equation}
is called a \emph{de Branges matrix} if $\det E_+ \not\equiv 0$, and $E_+^{-1}E_-$ belongs to $\mc{S}_{in}^{n \times n}$.
\end{definition}

\noindent An $n \times n$-matrix-valued entire function $E$ clearly belongs to HB$_n$ if and only if the restriction of $[E^{\#},E]$ to $\mathbb{C}_+$ is a de Branges matrix.
\medskip

\noindent Let us turn to the connection between de Branges matrices and meromorphic Herglotz functions. In analogy to \Cref{classichb}, let
\begin{align*}
A:=\frac{E_++E_-}{2}, \quad \quad \quad \quad B:=\frac{E_+-E_-}{2i}.
\end{align*}
Then, for $s:=E_+^{-1}E_-$,
\begin{align*}
E_+ \big[I-ss^* \big] E_+^* = E_+E_+^*-E_-E_-^*=2i \big(BA^*-AB^* \big).
\end{align*}

\noindent Under the assumption that $B(z)$ is invertible for every $z \in \bb{C}_+$, the function $Q:=B^{-1}A$ has positive semi-definite imaginary part at every $z \in \mathbb{C}_+$:
\begin{equation*}
\IM Q=\frac{1}{4} \big[B^{-1}E_+\big]\big[I-ss^* \big]\big[B^{-1}E_+\big]^* \geq 0.
\end{equation*}
If, in addition, $\mc{E}$ admits a meromorphic continuation to some open domain $\Omega \supset \bb{C}_+ \cup \bb{R}$, we can continue $Q$ to $\Omega$ as well (also denoting the continuation by $Q$). Then $\IM Q(x)$ vanishes for almost every real $x$, hence $Q(z)=Q^{\#}(z)$ whenever $z,\overline{z} \in \Omega$ and neither is a pole of $Q$. We can thus continue $Q$ to a $\#$-real, hence meromorphic Herglotz function defined on $\mathbb{C}$. In what follows, we will show that invertibility of $B$ is not necessary in order to define a meromorphic Herglotz function $Q$ in a meaningful way. However, we will require that $s=E_+^{-1}E_-$ admits a meromorphic continuation to an open and connected $\Omega \supset \bb{C}_+ \cup \bb{R}$. 
\newline
\noindent Because invertibility of $B$ is not given, we will work with the Moore-Penrose inverse $T^+$ of a matrix $T$. It is defined by the properties
\begin{itemize}
\item[$\bullet$] $TT^+T=T$,
\item[$\bullet$] $T^+TT^+=T^+$,
\item[$\bullet$] $TT^+$ and $T^+T$ are both Hermitian.
\end{itemize}
In fact, every matrix has a unique Moore-Penrose inverse. It is readily checked that $(T^+)^+=T$, and that $(\lambda T)^+=\frac{1}{\lambda}T^+$ for $\lambda \in \bb{C} \setminus \{0\}$. Additionally, $\ker T =\ker T^+T$ because of
\begin{equation}
\ker T \subseteq \ker T^+T \subseteq \ker TT^+T=\ker T.
\end{equation}
The Moore-Penrose inverse of a holomorphic function $F$ defined on a non-empty open and connected set $\Omega$ is holomorphic if and only if $\ran F$ and $\ker F$ are both constant on $\Omega$ (see \cite[Proposition 8.4]{fkms12}, \cite[Remark 3.8]{leiterer-rodman}).
\medskip

\noindent The next theorem is taken from \cite[Theorems 2.1 and 2.2]{fritzsche-k-m} and sheds light on the Moore-Penrose inverse of a Herglotz function.

\begin{theo}
\label{Herg_Moore_Penrose}
Let $Q$ be an $n \times n$-matrix-valued Herglotz function. Then $\ker Q(z)$ and $\ran Q(z)$ are independent of $z \in \mathbb{C} \setminus \mathbb{R}$. Take an orthonormal basis $(u_1,...,u_p)$ of $\ran Q(i)$ and let $U:=[u_1,u_2,...,u_p] \in \mathbb{C}^{n \times p}$ be the matrix whose columns are given by $u_1,...,u_p$. Then $q:=U^*QU$ is a $p \times p$-matrix-valued Herglotz function, and $q(z)$ is invertible for every $z \in \bb{C} \setminus \bb{R}$. Moreover, the following assertions hold:
\begin{itemize}
\item[$(i)$] $Q=UqU^*$;
\item[$(ii)$] $-Q^+=U[-q^{-1}]U^*$ is Herglotz;
\item[$(iii)$] $U^*U=I$ and $UU^*=Q(z)Q(z)^+=Q(z)^+Q(z)$ is the orthogonal projection onto $\ran Q(z)$ (which is independent of $z \in \mathbb{C} \setminus \mathbb{R}$).
\end{itemize}
\end{theo}

\begin{lemma}
\label{Herg_adjoint_kernel}
Let $M \in \bb{C}^{n \times n}$ such that $\IM M \geq 0$. Then $\ker M=\ker M^*$ and $\ran M=\ran M^*$.
\end{lemma}

\noindent A proof of \Cref{Herg_adjoint_kernel} is
given, e. g., in \cite[Definition 4.16 and Lemma A.32]{fkm21}. \\

\noindent For $s \in \mc{S}_{in}^{n \times n}$, introduce the function $R:=s-I$. Then, for $z \in \mathbb{C}_+$,
\begin{equation*}
0 \leq I-s(z)s(z)^*=-R(z)-R(z)^*-R(z)R(z)^*
\end{equation*}
and
\begin{equation}
\label{RisHerg}
2\IM (-iR(z))=-2\RE R(z)=-R(z)-R(z)^* \geq R(z)R(z)^* \geq 0. 
\end{equation}
Hence, the continuation of $-iR$ to a $\#$-real function on $\mathbb{C} \setminus \mathbb{R}$ is a Herglotz function, but in general we cannot continue it to a meromorphic Herglotz function (that is meromorphic on $\bb{C}$ and $\#$-real).
Observe that
\begin{equation}
B=-\frac{E_+R}{2i}, \quad \quad \quad \quad A=\frac{1}{2}E_+ \big(R+2I \big).
\end{equation}
If $B(z)$ is invertible for every $z \in \mathbb{C}_+$, it follows that $Q=B^{-1}A=-iR^{-1}\big(R+2I\big)$. If $B$ is not invertible, we define $Q:=-iR^+(R+2I)$ as an analog of $B^{-1}A$.
\medskip

\noindent The following theorem is reminiscent of Theorem 5.73 in \cite{arovdym}, in the sense that a connection between de Branges matrices and Herglotz functions is established. However, in \cite{arovdym} the associated Herglotz function has an absolutely continuous measure in its integral representation, while we are interested in the function $Q:=-iR^+(R+2I)$ which will extend to a meromorphic Herglotz function.

\begin{theo}
\label{hb_connection}
Let $\mc{E}=[E_-,E_+]: \bb{C}_+ \to \bb{C}^{n \times 2n}$ be meromorphic, such that $\det E_+ \not\equiv 0$, and $s=E_+^{-1}E_-$ admits a meromorphic continuation $\hat{s}$ to some open and connected $\Omega \supseteq \bb{C}_+ \cup \bb{R}$. Let $R=s-I$ and $Q=-iR^+(R+2I)$. Then the following assertions are equivalent:
\begin{itemize}
\item[$(i)$] $\mc{E}$ is a de Branges matrix;
\item[$(ii)$] $Q$ admits a continuation to a meromorphic Herglotz function on $\mathbb{C}$.
\end{itemize}
\end{theo}
\begin{proof} We prove first that $(ii)$ follows from $(i)$.
The calculation in (\ref{RisHerg}) shows that the $\#$-real continuation to $\mathbb{C} \setminus \mathbb{R}$ of $-iR$ is a Herglotz function. W.l.o.g., $R \not\equiv 0$. In view of \Cref{Herg_Moore_Penrose}, we can write $R(z)=Ur(z)U^*$ for all $z \in \mathbb{C}_+$, where $U \in \mathbb{C}^{n \times p}$ fulfills $U^*U=I$ and where $r : \mathbb{C}_+ \to \mathbb{C}^{p \times p}$ given by $r(z)=U^*R(z)U$ is such that $\IM (-ir(z)) \geq 0$ and $\det (-ir(z)) \neq 0$ for all $z \in \mathbb{C}_+$. Letting $\hat{R}:=\hat{s}-I$, we define a meromorphic continuation $\hat{r}$ of $r$ to $\Omega$, given by $\hat{r}(z)=U^*\hat{R}(z)U$ for all $z \in \Omega$. Define $T(z):=U \hat{r}(z)^{-1} U^*$ for $z \in \Omega$, which for $z \in \mathbb{C}_+$ coincides with $R(z)^+$. Then $\hat{R}T=UU^*$ and $T\hat{R}=UU^*$ are constant and equal to a Hermitian matrix. Since $T(z)\hat{R}(z)T(z)=T(z)$ for all $z \in \mathbb{C}_+$, we have $T\hat{R}T=T$ also on $\Omega$. Similarly, $\hat{R}T\hat{R}=\hat{R}$ on $\Omega$, showing that even $T(z)=\hat{R}(z)^+$ for all $z \in \Omega$. We see that $\hat{R}^+$ is meromorphic on $\Omega$. Hence, $Q$ admits a meromorphic continuation $\hat{Q}$ defined on $\Omega$. In addition, 
\begin{align}
\label{IM_hat_Q}
\IM \hat{Q} =-\hat{R}^+\hat{R}-\hat{R}^+-(\hat{R}^+)^*=\hat{R}^+ \big[-\hat{R}\hat{R}^*-\hat{R}^*-\hat{R} \big](\hat{R}^+)^* =\hat{R}^+ \big[I-\hat{s}\hat{s}^* \big](\hat{R}^+)^*.
\end{align}
In particular, $\IM \hat{Q}(z) \geq 0$ for all $z \in \mathbb{C}_+$ and $\IM \hat{Q} = 0$ a. e. on $\mathbb{R}$.  Consequently, $\hat{Q}$ can be continued to a meromorphic Herglotz function, and thus $(ii)$ holds true. 
\medskip

\noindent For the proof of the other implication, we first observe that $R$ is meromorphic on $\mathbb{C}_+$. If $\mc{P}(R)$ is the set of poles of $R$, then for every $z \in \mathbb{C}_+ \setminus \mc{P}(R)$, the matrix $R(z)^+R(z)$ is an orthogonal projection and, in particular, positive
semi-definite. For all $z \in \mathbb{C}_+ \setminus \mc{P}(R)$, we have
\begin{equation}
\label{R-Q-eq}
-iR(z)=-2\Big[Q(z)+iR(z)^+R(z)\Big]^+
\end{equation}
and that the matrix $Q(z)+ iR(z)^+R(z)$ has positive semi-definite $*$-imaginary part. Thus, from
(\ref{R-Q-eq}) and \Cref{Herg_Moore_Penrose} we get $\IM (-iR(z)) \geq 0$ for all $z \in \mathbb{C}_+ \setminus \mc{P}(R)$. This shows that $R$ has no poles in $\mathbb{C}_+$ and that the $\#$-real continuation of $-iR$ to $\mathbb{C} \setminus \mathbb{R}$ is a Herglotz function. \\
\noindent Consequently, because of \Cref{Herg_Moore_Penrose} and \Cref{Herg_adjoint_kernel}, $P := R(i)^+ R(i)$ is an orthogonal projection with $\ran P = \ran R(z) = \ran R(z)^*$
and $\ker P = \ker R(z) = \ker R(z)^*$ for all $z \in \mathbb{C}_+$. In view of $I - ss^*$ = $-RR^* - R^* - R$, on $\mathbb{C}_+$ we have $[I - ss^*](I - P ) = 0$ and $(I - P )[I - ss^* ] = 0$, hence
\begin{align}
I-ss^* &=P(I-ss^*)P =RR^+(I-ss^*)RR^+ \nonumber\\
\label{I_kernel_shrink}
&=(R^+)^* \big[R^*(I-ss^*)R \big] R^+. 
\end{align}
Taking into account
\begin{align}
\label{R-s-correspondence}
R^*(I-ss^*)R = R^*R \big[-R^+R-R^+-(R^+)^* \big] R^*R = R^*R (\IM Q) R^*R
\end{align}
then $\IM Q(z) \geq 0$ for $z \in \mathbb{C}_+$ implies $I-s(z)s(z)^* \geq 0$ for all $z \in \mathbb{C}_+$. If $R$ is replaced by $\hat{R}:=\hat{s}-I$ and $Q$ is replaced by its meromorphic continuation $\hat{Q}$ to $\Omega$, then, by continuity, (\ref{I_kernel_shrink}) and (\ref{R-s-correspondence}) also hold at every real $x$ that is not a pole of $\hat{Q}$ or $\hat{R}$. Since $\hat{Q}(z)=\hat{Q}^{\#}(z)$ whenever $z,\overline{z} \in \Omega$, we have $\hat{Q}(x)=\hat{Q}(x)^*$ and hence $\IM \hat{Q}(x) = 0$ for every such $x$. This shows that $I-\hat{s}\hat{s}^*=0$ a.e. on $\mathbb{R}$.
\end{proof}

\noindent It is also possible to start with a meromorphic Herglotz function $Q$ and parametrize the meromorphic functions $R$ such that $Q|_{\mathbb{C}_+}=-iR^+(R+2I)$. If such an $R$ is given, then from equation (\ref{R-Q-eq}) we see that $-iR$ can be continued to a Herglotz function, so $\ran R=\ran R^+ \supseteq \ran Q$. In fact, there is a one-to-one correspondence between such functions $R$ and linear subspaces containing $\ran Q$.

\begin{lemma}
Let $Q$ be an $n \times n$-matrix-valued Herglotz function and let $\mc{L}$ be a linear subspace of $\bb{C}^n$ that contains $\ran Q$. Then there is a unique holomorphic function $R$ on $\bb{C}_+$ such that $Q(z)=-iR(z)^+(R(z)+2I)$ and $\ran R(z)=\mc{L}$ for all $z \in \bb{C}_+$.
\end{lemma}
\begin{proof}
Let $P_{\mc{L}}$ be the orthogonal projection onto $\mc{L}$. Then a function $R$ satisfying the above conditions can be defined by
\begin{equation}
\label{R_proj_def}
R(z):=-2i\big(Q(z)+iP_{\mc{L}}\big)^+, \quad z \in \mathbb{C}_+.
\end{equation}
Indeed, $R$ is meromorphic because $(Q+iP_{\mc{L}})|_{\mathbb{C}_+}$ can be continued to a Herglotz function. By \Cref{Herg_Moore_Penrose}, $-iR$ admits continuation to a Herglotz  function, too. The theorem also states that $R^+R=RR^+$ is constant, and $\ker R(z)^+R(z)=\ker R(z)^+=\ker \big(Q(z)+iP_{\mc{L}}\big)$ for all $z \in \mathbb{C}_+$. We know from \Cref{Herg_adjoint_kernel} that 
\begin{align*}
&\ker \big(Q(z)+iP_{\mc{L}}\big) \\
&=\big(\ker \RE Q(z) \big) \cap \Big(\ker \big(\IM Q(z)+P_{\mc{L}}\big)\Big) \\
&=\ker (\RE Q(z)) \cap \ker (\IM Q(z)) \cap \ker P_{\mc{L}}  \\
&=\ker Q(z) \cap \mc{L}^{\perp}=\mc{L}^{\perp}, \quad z \in \mathbb{C}_+.
\end{align*}
Clearly, $R^+R$ is an orthogonal projection, and it follows that $R^+R=P_{\mc{L}}$. Putting $R^+R$ in place of $P_{\mc{L}}$ in (\ref{R_proj_def}) implies $Q|_{\mathbb{C}_+}=-iR^+(R+2I)$. Uniqueness of $R$ follows in the same way.
\end{proof}

\noindent In the setting of the above lemma, if $Q$ is even a meromorphic Herglotz function, $R$ admits a meromorphic continuation to $\bb{C}$: In view of \Cref{Herg_Moore_Penrose}, pick an orthonormal basis $(u_1,...,u_p)$ of $\mc{L}=\ran (Q(i)+iP_{\mc{L}})$ and let $U=[u_1,u_2,...,u_p]$. Then $\tilde{q}:=U^*(Q+iP_{\mc{L}})U$ is a $p \times p$-matrix-valued meromorphic function on $\bb{C}$, and $\tilde{q}(z)$ is invertible for every $z \in \bb{C}_+$. It follows that $\hat{R}:=-2iU\tilde{q}^{-1}U^*$ is a meromorphic continuation of $R$ to $\bb{C}$. Looking at the statement and the proof of \Cref{hb_connection}, we get as a corollary that an $s \in \mc{S}_{in}^{n \times n}$ admits a meromorphic continuation to all of $\bb{C}$ if and only if it admits a meromorphic continuation to some open and connected $\Omega \supset \bb{C}_+ \cup \bb{R}$.
\medskip

\noindent When combining \Cref{sufficiency} with \Cref{hb_connection}, one obtains the following statement, which can be viewed as a matrix-valued version of the Hermite-Biehler Theorem.

\begin{theo}
\label{final-matrix-hb}
Let $\mc{E}=[E_-,E_+]: \bb{C}_+ \to \bb{C}^{n \times 2n}$ be meromorphic, such that $\det E_+ \not\equiv 0$, and $s=E_+^{-1}E_-$ admits a meromorphic continuation $\hat{s}$ to some open and connected $\Omega \supseteq \bb{C}_+ \cup \bb{R}$. Let $R=s-I$, $\hat{R}=\hat{s}-I$, and $\hat{Q}=-i\hat{R}^+(\hat{R}+2I)$.
Then $\mc{E}$ is a de Branges matrix if and only if
\begin{itemize}
\item[$(i)$] $R$ is of bounded type;
\item[$(ii)$] Except for poles, $\ker R(z)$ is independent of $z \in \bb{C}_+$, and $R^+$ has no poles in $\mathbb{C}_+$;
\item[$(iii)$] $I-\hat{s}(x)\hat{s}(x)^*=0$ for a.e. $x \in \mathbb{R}$;
\item[$(iv)$] $\limsup_{\eta \to +\infty} \frac{\|R(i\eta)^+\|}{\eta}<+\infty$.
\item[$(v)$] For every $m \in \{1,...,n\}$ and $\mathbf{i} \in I_{n,m}$, the function $f_{\mathbf{i}}:=\det \hat{Q}_{\mathbf{i}}$ satisfies the $m$-interlacing condition. If $f_{\mathbf{i}}$ has a pole of multiplicity $m$, there is $z^* \in \mathbb{R}$ such that $\lim_{z \to z^*} [(z^*-z)^mf_{\mathbf{i}}(z)]>0$. The limit $\lim_{\eta \to +\infty} \frac{f_{\mathbf{i}}(i\eta)}{(i\eta)^m}$, which exists by \Cref{colour_func}, is nonnegative.
\end{itemize}
\end{theo}
\begin{proof}
If $\mc{E}$ is a de Branges matrix, the first four properties follow from the previous results, and the last one from \Cref{sufficiency}. For the reverse implication, we need to check that $\hat{Q}$ extends to a meromorphic Herglotz function, then apply \Cref{hb_connection}. Note that by (\ref{IM_hat_Q}) and by $(iii)$, we have $\IM \hat{Q}(x)=0$ for almost every $x \in \mathbb{R}$. We may thus assume that $\hat{Q}$ is a $\#$-real meromorphic function on all of $\mathbb{C}$. Observe that $R(z)^+R(z)$ is the orthogonal projection onto $\ran R(z)^+=(\ker R(z))^{\perp}$, which by assumption is independent of $z \in \mathbb{C}_+$. Therefore, $R^+R$ is constant. Since we also assumed that $R^+$ has no poles in $\mathbb{C}_+$, the same holds for $\hat{Q}$. Hence, \Cref{sufficiency} is applicable and shows that $\hat{Q}$ is Herglotz.
\end{proof}

\subsection*{Acknowledgements.}
\noindent The author wants to thank Roman Romanov and Harald Woracek for their support and guidance.

\section{Declarations}
\noindent \textbf{Funding} This work was supported by the Austrian Science Fund [grant number I-4600]. \\
\textbf{Conflicts of interest:} \\
The author has no relevant interests to disclose.

\end{document}